\documentclass{birkjour}
\usepackage{amsxtra}
\usepackage{amsthm,mathrsfs}
\usepackage{amsmath,amsthm,amssymb,amsrefs,bbm,color,enumitem,esint,esvect,float,graphicx,mathrsfs}
\begin{document}
\begin{abstract}
We resolve the question of the boundedness of the composition of dyadic paraproducts, first posed by Pott, Reguera, Sawyer, and Wick in~\cite{PotCarSawWic}, by providing necessary and sufficient conditions for their boundedness.
\end{abstract}
 \newtheorem{theorem}{Theorem}[section]
 \newtheorem{corrolary}[theorem]{Corollary}
 \newtheorem{lemma}[theorem]{Lemma}
 \newtheorem{proposition}[theorem]{Proposition}
 \theoremstyle{definition}
 \newtheorem{definition}[theorem]{Definition}
 \theoremstyle{remark}
 \newtheorem{remark}[theorem]{Remark}
 \newtheorem*{example}{Example}
 \numberwithin{equation}{section}

\title{Boundedness of a composition of paraproducts}
\author{Ana Čolović}
\address{
Mathematics and Science Building\\
University of Missouri\\
810 Rollins St\\
Columbia, MO USA\\
65201}
\email{acg7y@umsystem.edu}

\maketitle

\section{Introduction}

Paraproduct operators are the building blocks of the multiplication operator $M_b:L^2(\mathbb{R})\to L^2(\mathbb{R}),$ for $b\in L^1_{\text{loc}, 
}$ and are bounded under a more general set of conditions than the multiplication operator (see~\cite{BenMalNai}). While the operator $M_b$ is bounded on $L^2(\mathbb{R})$ if and only if $b\in L^{\infty},$ the paraproduct operator $\Pi_b,$ is bounded on $L^2(\mathbb{R}),$ if and only if $b\in BMO,$ the space of bounded mean oscillation. Given their role in decomposing the multiplication operator $M_b$, they arise as an essential component of understanding commutators of Calder\'{o}n-Zygmund operators $T,$ and the multiplication operator $M_b,$ on $L^p$ spaces, where $1<p<\infty,$ usually denoted by $[M_b,T],$ or simply as $[b,T]$, first in the work of Calder\'{o}n (in~\cite{Cal}), and then later in more general weighted settings (see~\cite{Chu}).

More broadly, dyadic paraproducts have been used to model similar operators in complex analysis, namely Toeplitz and Hankel operators (see~\cite{PotCarSawWic}). 

For a fixed symbol function $\varphi,$ the Toeplitz operator $T_{\varphi}:H^2(\mathbb{T})\to H^2(\mathbb{T})$ is the operator defined  by $T_{\varphi}f :=\mathbb{P}_{+}(\varphi f),$ where $\mathbb{P}_{+}$ is the projection from $L^2(\mathbb{T})$ to the Hardy space $H^2(\mathbb{T}),$ the space of all functions $f\in L^2(\mathbb{T})$ with vanishing negative Fourier coefficients. Since Toeplitz operators are defined in terms of a projection of a product of two functions, one can think of them as analogous to paraproduct operators in the analytic setting. They are bounded if and only if the following condition holds on the symbol:
\[
\|T_{\varphi}\|_{H^2(\mathbb{T})}=\sup_{z\in \mathbb{D}}|\varphi(z)|<\infty,
\]
where $\varphi(z)$ denotes the harmonic extension of the function $\varphi$ evaluated at a point $z\in \mathbb{D}.$

 The study of compositions of unbounded Toeplitz operators was started by Sarason (in~\cite{Sar}) who conjectured that the composition  $T_{\varphi} \circ T_{\overline{\psi}}$ is bounded if and only if
$\sup_{z\in \mathbb{D}}|\varphi|^2(z) |\psi|^2(z)< \infty,$ where $|\varphi|(z)$ and $|\psi|(z)$ are harmonic extensions of $|\varphi|$ and $|\psi|$ evaluated at $z\in \mathbb{D}.$  The conjecture combined the two conditions for the boundedness of each of the Toeplitz operators, and posited that their composition will remain bounded under a simple joint condition. The conjecture was disproved by Nazarov, who in~\cite{Naz} found a counterexample. One could still ask when will a composition of two Toeplitz operators be bounded. Cruz-Uribe, in~\cite{Cru}, showed that the boundedness of a composition of Toeplitz operators is equivalent to the two-weight boundedness of the Hilbert transform on $\mathbb{R}$. In particular, the composition $T_{\varphi}\circ T_{\overline{\psi}}$ is bounded on $L^2(\mathbb{T})$ if and only if the Hilbert transform is bounded as an operator between $L^2(\mathbb{T}, |\psi|^{-2})$ and $L^2(\mathbb{T}, |\phi|^2).$ 
A longstanding question in harmonic analysis, the problem of two-weight boundedness of the Hilbert transform was solved in two papers, by Lacey (in~\cite{Lac}) and Lacey, Sawyer, Shen, Uriarte-Tuero (in~\cite{LacSawSheUri}). Given the equivalence of the boundedness of the Hilbert transform to that of boundedness of a composition of Toeplitz operators, theorems in~\cite{Lac} and~\cite{LacSawSheUri} gave a set of conditions under which the composition of Toeplitz operators is bounded on $L^2(\mathbb{T}).$ 

 Given the analogous role in decomposing products of two functions, dyadic paraproducts can be thought of as discrete analogues of the Toeplitz operators. Because of a similar role they play, and the connection of compositions of Toeplitz operators to important questions in harmonic analysis, in~\cite{PotCarSawWic}, Pott, Reguera, Sawyer and Wick began the study of the the composition of dyadic paraproduct operators. 

For an interval $I\in \mathcal{D},$ where $\mathcal{D}$ is the standard grid of dyadic intervals, we define 
\[
h_I^{0} \equiv h_I \equiv \frac{1}{\sqrt{|I|}} (\mathbf{1}_{I^{+}}-\mathbf{1}_{I^{-}}) \; \text{ and } \; h^{1}_{I} \equiv \frac{\mathbf{1}_{I}}{|I|}, \; I\in \mathcal{D}.
\]

 Given a sequence of complex numbers $b=\{b_I\}_{I \in \mathcal{D}},$ one can define operators known as dyadic paraproducts, as in~\cite{PotCarSawWic}.

\begin{definition} Let $b=\{b_{I}\}_{I\in \mathcal{D}}$, and $(\alpha, \beta) \in \{0,1\} \times \{0,1\}.$ The dyadic paraproduct acting on a function $f$ is given by: 
\[
P^{(\alpha,\beta)}_{b} f \equiv \sum_{I\in \mathcal{D}} b_{I} \langle f, h_{I}^{\beta}\rangle_{L^2(\mathbb{R})} h_{I}^{\alpha}.
\] We refer to $(\alpha, \beta)$ as the index of the paraproduct. The sequence $b=\{b_I\}_{I\in \mathcal{D}}$ is referred to as a symbol of the operator $P_b^{(\alpha, \beta)}.$

\end{definition}

The operator $P_b^{(0,1)}$ is often denoted by $\Pi_b.$ Its adjoint in $L^2(\mathbb{R})$ is the operator $P^{(1,0)}_b = \Pi_b^*.$ 

In~\cite{PotCarSawWic}, the authors considered the boundedness of compositions of certain types of paraproduct operators $P_b^{(\alpha, \beta)}\circ P_{d}^{(\gamma, \delta)}$ on $L^2(\mathbb{R}).$ In~\cite{PotCarSawWic}, the $P_b^{(\alpha, \beta)}\circ P_{d}^{(\gamma, \delta)}$ is referred to as $P_{b,d}^{(\alpha, \beta, \gamma, \delta)},$ and the index $(\alpha, \beta, \gamma, \delta)$ is referred to as the type of the composition of two paraproducts. We note that, for a sequence $b$ an adjoint of the operator $P_{b}^{(\alpha,\beta)},$ is given by $P^{(\beta,\alpha)}_{\overline{b}}.$
Pott, Reguera, Sawyer and Wick (in~\cite{PotCarSawWic}) characterized the conditions for boundedness of type $(\alpha,0,0,\beta)$ and $(0,\alpha,\beta, 0)$ operators.
In general, the adjoint of the paraproduct operator $P_{b}^{(\alpha, \beta)}$ is the operator $P_{\overline{b}}^{(\beta, \alpha)}$, where $(\alpha,\beta)\in \{0,1\}^2,$ so once the boundedness of operators studied by~\cite{PotCarSawWic} was established, the remaining operator to characterize was the operator of type $(0,1,0,1),$ which arises as a composition of operators $\Pi_b$ and $\Pi_d,$ for symbols $b$ and $d.$

In this paper, we resolve the remaining case, giving conditions for the boundedness of the operator $\Pi_b\Pi_d.$ We prove the following theorem.

\begin{theorem} \label{T:0101}
Let $b=\{b_I\}_{I\in \mathcal{D}},$ $d=\{d_I\}_{I\in \mathcal{D}}.$
For $I\in \mathcal{D},$ let \[
\nu(I):=\sum_{J\subset I} |b_J|^2, \: \text{ and }\: \widetilde{\nu}(I):= \sqrt{\frac{\nu(I^{+})\nu(I^{-})}{\nu(I^+)+\nu(I^{-})}}.
\] 
For $I,J\in \mathcal{D}, J \subsetneq I,$ let $\delta^{(I,J)}=\begin{cases} 
      1 & J\subset I^{+} \\
      -1 & J\subset I^- \\
      0  & I\cap J = \emptyset.
   \end{cases}$

    The operator $P_{b,d}^{(0,1,0,1)}$ is bounded on $L^2(\mathbb{R})$ if and only if the following conditions hold:
    \begin{enumerate}
     \item \begin{equation} \label{E:Condition-A}\sup_{I\in \mathcal{D}} \left|\frac{1}{\sqrt{|I|}} \sum_{J\subset I^{\pm}} \frac{d_J}{\sqrt{|J|}}\frac{\widetilde{\nu}(I)}{\nu(I^{\pm})}(-\nu(J^{+})+\nu(J^{-})) \right|:=A< \infty\end{equation}
         \item \begin{align} \label{E:Condition-B}&\sup_{I\in\mathcal{D}} \left(\frac{1}{|I|}\sum_{J\subsetneq I}\left(\widetilde{\nu}(J)\right)^2\left(\sum_{K\subset J^{\pm}} \frac{d_K}{\sqrt{|K|}}\frac{\delta^{(J,K)}}{\nu(J^{\pm})}\left(-\nu(K^+)+\nu(K^-)\right) + \frac{d_J}{\sqrt{|J|}}\right)^2\right)^{\frac{1}{2}} \\ &:= B < \infty.\end{align} 
        \item \begin{equation} \label{E:Condition-C}\sup_{I\in\mathcal{D}}\left(\frac{1}{\nu(I)}\sum_{J\subset I}\frac{1}{|J|} \left(\sum_{K\subset J^{\pm}} \delta^{(J,K)} \frac{d_K}{\sqrt{|K|}} (-\nu(K^+)+\nu(K^-))\right)^2\right)^{\frac{1}{2}}:= C< \infty \end{equation}
        
    \end{enumerate}
    The norm of $\|\Pi_b\Pi_d:L^2(\mathbb{R})\to L^2(\mathbb{R})\|\simeq A+B+C.$
    Moreover, $(A + B + C) \lesssim \|b\|_{BMO} \|d\|_{BMO}.$
\end{theorem}
We show that $P_{b,d}^{(0,1,0,1)}$ has the same boundedness condition as a certain operator $T_{b,d}^{(0,1,0,1 )},$ defined on the upper half-plane. Conditions~\eqref{E:Condition-A}, \eqref{E:Condition-B}, and \eqref{E:Condition-C} arise when the boundedness of the operator $T_{b,d}^{(0,1,0,1)}$ is tested on appropriate functions defined on the upper half-plane. We use an argument similar to that presented in~\cite{PotCarSawWic}, which modifies the approach of~\cite{TreVol} to well-localized operators. In the following section, we introduce the relevant spaces on the upper half-plane.  In Section \ref{S:Proof-0101}, we provide the proof of Theorem \ref{T:0101}, along with the proof that the bounds constructed are at least as good as the $BMO$ bounds.

\subsection{Acknowledgments} The author would like to thank Brett Wick for many useful discussions concerning this paper. 

\section{Background on the Upper Half-Plane} 

To argue that the operator $\Pi_b \, \Pi_d$ is bounded, we transplant it to the upper half-plane. The objects that appear in this section can also be found in~\cite{PotCarSawWic}, but we introduce them for ease of reference and clarity. 

    We use $\mathcal{H}$ to denote the upper half-plane. To each interval $I$ in the dyadic grid $\mathcal{D},$ on the real line we associate a Carleson tile on the upper half-plane defined by:
  $T(I):= I \times \left[\frac{|I|}{2}, |I|\right].$ A Carleson cube corresponding to the interval $I$ is defined by: $
    Q(I):= I \times \left[0,|I|\right].$
    Carleson tiles are disjoint and their union is the entire upper half-plane, i.e. $\mathcal{H}:=\bigcup_{I\in \mathcal{D}}T(I).$ For each $I \in \mathcal{D},$ $Q(I)=\bigcup_{J\subset I}T(J).$

We let $L^2(\mathcal{H})$ denote the space of square integrable functions on the upper half-plane. For a function $\sigma \geq 0,$ we let $L^2(\mathcal{H},\sigma)$ denote the space of functions $f$ on the upper half-plane,
such that 
$\|f\|^2_{L^2(\mathcal{H},\sigma)}=\int_{\mathcal{H}}|f(z)|^2 \sigma(z) dA(z) < \infty.$

The subspace $L^2_c(\mathcal{H})$ denotes the space of functions in the upper half-plane constant on the Carleson tiles or equivalently functions $f:\mathcal{D}\to \mathbb{C}$,
\[
f=\sum_{I\in \mathcal{D}}f_I \mathbf{1}_{T(I)}
\]
such that 
\[
\|f\|^2_{L^2_c(\mathcal{H})}:= \frac{1}{2}\sum_{I\in \mathcal{D}} |f(I)|^2 |I|^2 < \infty.
\]

The indicator functions of the Carleson tiles are orthogonal in $L^2_c(\mathcal{H}).$
For a function $f\in L^2(\mathcal{H}),$ as in~\cite{PotCarSawWic}, 
we use $\widetilde{f}$ to denote $\frac{f}{\|f\|_{L^2(\mathcal{H})}}.$ Note that  $\{\widetilde{\mathbf{1}}_{T(I)}\}_{I\in \mathcal{D}}$ is then an orthonormal basis of $L^2_c(\mathcal{H}).$ We observe that for $I\in \mathcal{D},$
\[
\widetilde{\mathbf{1}}_{T(I)}= \frac{\mathbf{1}_{T(I)}}{\|\mathbf{1}_{T(I)}\|_{L^2(\mathcal{H})}}=\frac{\mathbf{1}_{T(I)}}{\sqrt{|T(I)|}}= \sqrt{2}\frac{\mathbf{1}_{T(I)}}{|I|}.
\]

Let $b=\{b_I\}_{I\in \mathcal{D}},$ be a sequence of complex numbers, and $\alpha \in \mathbb{R}.$ We define the operator $\mathcal{M}_b^{\alpha}$ on $L^2_c(\mathcal{H})$ by its action on $\widetilde{\mathbf{1}}_{T(I)}$ by
\[
\mathcal{M}_{b}^{\alpha}\widetilde{\mathbf{1}}_{T(I)}= b_I |I|^{\alpha} \widetilde{\mathbf{1}}_{T(I)}.
\]

For a non-negative function $w$ on the dyadic grid $\mathcal{D},$ we let $l^2(w)$ denote the space of sequences $\{b_I\}_{I\in \mathcal{D}},$ such that 
\[
\sum_{I\in \mathcal{D}}|b_I|^2 w(I) < \infty.
\]
These spaces of functions allow for us to transplant the problem of the boundedness of the operator $P^{(0,1,0,1)}_{b,d},$ with symbols $b,d,$ into the question of boundedness of an operator on $L^2(\mathcal{H}).$ 

\section{Boundedness of Type (0,1,0,1) Composition} \label{S:Proof-0101}

Using the transplantation method to the upper half-plane developed by Pott,Reguera, Sawyer and Wick in~\cite{PotCarSawWic}, we characterize the boundedness for the operator $P_b^{(0,1)}\circ P_d^{(0,1)}.$ We observe that the operator $P_b^{(0,1)}\circ P_d^{(0,1)}$ can be written as a composition of two operators characterized by~\cite{PotCarSawWic}. This allows for the transplantation to the upper half-plane, and characterizing the boundedness of the operator $P_b^{(0,1)}\circ P_d^{(0,1)}.$
\begin{proof}{Proof of Theorem \ref{T:0101}} We consider the Gram matrix $[G_{I,J}]_{I,J\in\mathcal{D}}$ of the operator $P_b^{(0,1)}\circ P_d ^{(0,1)}$ relative to the Haar basis.
We note that 
\begin{align*}
G_{I,J} &= \langle P^{(0,1)}_b\circ P_d^{(0,1)}h_J, h_I\rangle_{L^2(\mathbb{R})} = \langle P_d^{(0,1)}h_J, P_b^{(1,0)}h_I\rangle_{L^2(\mathbb{R})} \\
& = \langle \sum_{K\subsetneq J} d_{K} h_J(K) h_K, b_I h^1_I\rangle_{L^2(\mathbb{R})} \\
&= \sum_{I\subsetneq K \subsetneq J} d_{K}\overline{b_I} h_J(K)h_{K}(I).
\end{align*}
In~\cite{PotCarSawWic}, the authors prove that the operator $P^{(0,1,0,0)}_{e,f}$ for symbols $e$ and $f,$ has the same Gram matrix with respect to the Haar basis as an operator on the upper half-plane. The identity operator $I$ on $L^2(\mathbb{R})$ can be written as 
$If=\sum_{I\in \mathcal{D}}\langle f, h_I\rangle h_I.$ Hence, the identity operator on $L^2(\mathbb{R})$ can be thought of as a dyadic martingale operator $P^{(0,0)}_{\epsilon},$ with a symbol $\epsilon =1. $ Therefore, we can factor $P_b^{(0,1)}\circ P_d^{(0,1)}$ as 
$P_b^{(0,1)}\circ P_{1}^{(0,0)}\circ P_d^{(0,1)}\circ P_1^{(0,0)}.$ 

In~\cite{PotCarSawWic}, the authors find an operator on the upper half-plane that has the same Gram matrix with respect to an appropriate basis as the Gram matrix of the operator $P_{e,f}^{(0,1,0,0)},$ with the basis of Haar wavelets. We use the transplantation of~\cite{PotCarSawWic} and the factorization of the operator $P_{b,d}^{(0,1,0,1)}$ into two operators of the $(0,1,0,0)$ type to transplant the operator $P_{b,d}^{(0,1,0,1)}$ onto the upper half-plane.

Let \[T^{(0,1,0,1)}_{b,d} := \mathcal{M}_{\overline{b}}^{-1}U\mathcal{M}_{\overline{d}}^{-\frac{1}{2}}U\mathcal{M}_{1}^\frac{1}{2},\]
where \[U=\sum_{J\in\mathcal{D}}\widetilde{\mathbf{1}}_{Q_{\pm}(J)}\otimes \widetilde{\mathbf{1}}_{T(J)},\]
and \[\mathbf{1}_{Q_{\pm}(J)}:= - \sum_{K\subset J_-}\mathbf{1}_{T(K)}+\sum_{K\subset J_+}\mathbf{1}_{T(K)}.\]

Note that for $J\in \mathcal{D},$
\begin{equation}\label{E:Norm-of-Q+-}
    \|\mathbf{1}_{Q_{\pm}(J)}\|_{L^2(\mathcal{H})}=\frac{|J|}{2}.
\end{equation}

Let $K\in \mathcal{D}.$ 
For $\alpha\in\mathbb{R},$ a symbol $a,$ and $K\in \mathcal{D}$, we note that 
\[
\mathcal{M}^{\alpha}_{a}\mathbf{1}_{Q_{\pm}(K)}=\sum_{L\subset K^{\pm}} \delta^{(K,L)} a_L |L|^{\alpha} \mathbf{1}_{T(L)}.
\]

The entries of the matrix $[G_{I,J}]_{I,J\in\mathcal{D}}$ of $T_{b,d}^{(0,1,0,1)}$ with respect to the basis $\{\widetilde{\mathbf{1}}_{T(I)}\}_{I\in \mathcal{D}}$ are
    \begin{align*}
    G_{I,J}&=\langle \mathcal{M}_{\overline{b}}^{-1}U\mathcal{M}_{\overline{d}}^{-\frac{1}{2}}U\mathcal{M}_{1}^\frac{1}{2} \widetilde{\mathbf{1}}_{T(J)} , \widetilde{\mathbf{1}}_{T(I)}\rangle_{L^2(\mathcal{H})} \\
    &= \langle \mathcal{M}_{\overline{b}}^{-1}U\mathcal{M}_{\overline{d}}^{-\frac{1}{2}}|J|^{\frac{1}{2}} U \widetilde{\mathbf{1}}_{T(J)} , \widetilde{\mathbf{1}}_{T(I)}\rangle_{L^2(\mathcal{H})}  \\
    &= \langle \mathcal{M}_{\overline{b}}^{-1}U \mathcal{M}_{\overline{d}}^{-\frac{1}{2}}|J|^{\frac{1}{2}} \widetilde{\mathbf{1}}_{Q_{\pm}(J)} , \widetilde{\mathbf{1}}_{T(I)}\rangle_{L^2(\mathcal{H})} \\
    &= \langle \mathcal{M}_{\overline{b}}^{-1} |J|^{\frac{1}{2}} \sum_{K\in \mathcal{D}} (\widetilde{\mathbf{1}}_{Q_{\pm}(K)} \otimes \widetilde{\mathbf{1}}_{T(K)}) \mathcal{M}_{\overline{d}}^{-\frac{1}{2}} \widetilde{\mathbf{1}}_{Q_{\pm}(J)} , \widetilde{\mathbf{1}}_{T(I)}\rangle_{L^2(\mathcal{H})} \\
    &= \langle \mathcal{M}_{\overline{b}}^{-1} |J|^{\frac{1}{2}} \sum_{K\in \mathcal{D}} \langle \widetilde{\mathbf{1}}_{T(K)},  \mathcal{M}_{\overline{d}}^{-\frac{1}{2}} \widetilde{\mathbf{1}}_{Q_{\pm}(J)} \rangle_{L^2(\mathcal{H})}  \widetilde{\mathbf{1}}_{Q_{\pm}(K)}, \widetilde{\mathbf{1}}_{T(I)}\rangle_{L^2(\mathcal{H})} \\
     &= \sum_{K\in \mathcal{D}} |J|^{\frac{1}{2}} \langle \mathcal{M}_{d}^{-\frac{1}{2}} \widetilde{\mathbf{1}}_{Q_{\pm}(J)}, \widetilde{\mathbf{1}}_{T(K)}   \rangle _{L^2(\mathcal{H})}  \langle \mathcal{M}_{\overline{b}}^{-1}  \widetilde{\mathbf{1}}_{Q_{\pm}(K)}, \widetilde{\mathbf{1}}_{T(I)}\rangle_{L^2(\mathcal{H})}
     \end{align*}
We make note of the following calculation, used in computing the Gram matrix coefficients in the equation~ \eqref{E:Gram-matrix-coefficients}:
\begin{align} \label{E:Inner-product-Q+-}
&\langle\mathcal{M}_{a}^{\alpha}\widetilde{\mathbf{1}}_{Q_{\pm}(K)}, \widetilde{\mathbf{1}}_{T(I)}\rangle_{L^2(\mathcal{H})} \notag\\&= \frac{1}{\|\mathbf{1}_{Q_{\pm}(K)}\|_{L^2(\mathcal{H})}\|\mathbf{1}_{T(I)}\|_{L^2(\mathcal{H})}}\langle \sum_{L\subset K^{\pm}} \delta^{(K,L)} a_L |L|^{\alpha} \mathbf{1}_{T(L)}, \mathbf{1}_{T(I)}\rangle_{L^2(\mathcal{H})} \notag \\
&=\frac{1}{\|\mathbf{1}_{Q_{\pm}(K)}\|_{L^2(\mathcal{H})}\|\mathbf{1}_{T(I)}\|_{L^2(\mathcal{H})}} \delta^{(K,I)} a_I |I|^{\alpha}|T(I)|\notag \\
&= \sqrt{2}a_I \delta^{(K,I)} \frac{|I|^{\alpha+1}}{|K|},
\end{align}
for $I,K\in \mathcal{D},$ $\alpha \in \mathbb{R},$ and a symbol $a$.
Using the equation~\eqref{E:Inner-product-Q+-}, we then have that
\begin{align} \label{E:Gram-matrix-coefficients}
     G_{I,J}&=\sum_{K\in \mathcal{D}} |J|^{\frac{1}{2}} \langle \mathcal{M}_{d}^{-\frac{1}{2}} \widetilde{\mathbf{1}}_{Q_{\pm}(J)}, \widetilde{\mathbf{1}}_{T(K)}   \rangle_{L^2(\mathcal{H})}   \langle \mathcal{M}_{\overline{b}}^{-1}  \widetilde{\mathbf{1}}_{Q_{\pm}(K)}, \widetilde{\mathbf{1}}_{T(I)}\rangle_{L^2(\mathcal{H})}\\
     &= 2 \sum_{\substack{K\subset J^{\pm} \notag \\ I \subset K^{\pm}}} \left(\delta^{(J,K)} d_K |J|^{\frac{1}{2}}|J|^{-1}|K|^{-\frac{1}{2}}\right) \left( \overline{b_I}\frac{|I|^{-1+1}}{|K|}\right) \notag \\
     &= 2 \sum_{\substack{K\subset J^{\pm} \\ I \subset K^{\pm}}} d_K \overline{b_I} \frac{\delta^{(J,K)}}{|J|^{\frac{1}{2}}} \frac{\delta^{(I,K)}}{|K|^{\frac{1}{2}}} \notag
    \end{align}
Thus, the coefficients of the matrix of the operator $T_{b,d}^{(0,1,0,1)}$ with respect to the basis $\{\widetilde{\mathbf{1}}_{T(I)}\}_{I\in \mathcal{D}}$ agree, up to the constant 2, with those of the operator $P_{b,d}^{(0,1,0,1)}$ with respect to the basis of Haar coefficients. Therefore,
\[
\|P^{(0,1,0,1)}_{b,d}\|_{L^2(\mathbb{R})} \approx \|T^{(0,1,0,1)}_{b,d}\|_{L^2(\mathcal{H})}. 
\]

Following the method of~\cite{PotCarSawWic}, we rephrase the boundedness of the operator \\ $T^{(0,1,0,1)}_{b,d}$ on $L^2(\mathcal{H})$ in terms of a boundedness of a related, but simpler to analyze operator on a weighted space. We can then characterize the boundedness of $T^{(0,1,0,1)}_{b,d}$ in terms of testing conditions on $L^2(\mathcal{H})$, leading us to the boundedness conditions for the operator $P_{b,d}^{(0,1,0,1)}.$

The operator $T^{(0,1,0,1)}_{b,d}$ is bounded if and only if
\begin{equation}\label{E:Upper-operator-ineq}
\|\mathcal{M}_{\overline{b}}^{-1}U\mathcal{M}_{\overline{d}}^{-\frac{1}{2}}U\mathcal{M}_{1}^\frac{1}{2}f\|_{L^2_c(\mathcal{H})} \lesssim \|f\|_{L^2_c(\mathcal{H})},  \: \text{for all } \: f\in L^{2}(\mathcal{H}).
\end{equation}
 
For $f \in L^2_c(\mathcal{H}),$ let $\mathcal{M}^{\frac{1}{2}}_{1}f=g.$ Then the equation~\eqref{E:Upper-operator-ineq} becomes:

\begin{equation}\label{E:Upper-operator-weighted}
    \|U\mathcal{M}_{\overline{d}}^{-\frac{1}{2}}Ug\|_{L^2_c(\mathcal{H,\nu})} \lesssim \|g\|_{L^2_c(\mathcal{H},w)},
\end{equation}
where
\begin{align*}
&\nu := \sum_{I\in \mathcal{D}} |b_I|^2|I|^{-2} \mathbf{1}_{T(I)}\\
& w:= \sum_{I\in \mathcal{D}} |I|^{-1} \mathbf{1}_{T(I)}.
\end{align*}
We let \[
\mu:= \sum_{I\in \mathcal{D}} |I| \mathbf{1}_{T(I)},
\]
and substitute $g=\mu \,h,$ into the equation~\eqref{E:Upper-operator-weighted} to get
\begin{equation} \label{E:Upper-final-form}
    \|U\mathcal{M}_{\overline{d}}^{-\frac{1}{2}}U(\mu\, h)\|_{L^2_c(\mathcal{H,\nu})} \lesssim \|h\|_{L^2_c(\mathcal{H},\mu)}.
\end{equation}
Thus, to characterize the boundedness of the operator $T_{b,d}^{(0,1,0,1)},$ one needs to find the conditions under which the equation~\eqref{E:Upper-final-form} holds.

\subsection{Necessity of the Conditions}
We now show that the conditions $A$, $B,$ and $C$ that appear in Theorem \ref{T:0101} arise from testing the operator  $U\mathcal{M}_{\overline{d}}^{-\frac{1}{2}}U(\mu \, \cdot),$ and its adjoint on appropriate functions in $L^2_c(\mathcal{H},\mu)$ and $L^2_c(\mathcal{H}, \nu).$ In other words, we prove that the conditions $A, B$ and $C$ are necessary for boundedness of the operator $UM_{\overline{d}}^{-\frac{1}{2}}U(\mu \,\cdot).$

To do so, we use the orthonormal bases for $L^2_c(\mathcal{H}, \mu)$ and $L^2_c(\mathcal{H}, \nu),$ introduced in~\cite{PotCarSawWic}. We follow the notation of~\cite{PotCarSawWic} and use $\nu(J)$ to denote $\nu(Q(J)).$ For $I\in \mathcal{D}, $  we let $\mu_I=\mu(T(I)).$ The corresponding orthonormal bases for $L^2(\mu)$ and $L^2(\nu)$ are:
\[
\{h_I^{\mu}\}_{I\in\mathcal{D}}, \quad\{H_J^{\nu}\}_{J\in\mathcal{D}},
\]
where 
\begin{equation} \label{E:weighted-bases}
h_{I}^{\mu}:= \frac{\widetilde{\mathbf{1}}_{T(I)}}{\sqrt{\mu_I}}, \: \text{ and } \: H_{J}^{\nu}:= \widetilde{\nu}(J) \left(-\frac{\mathbf{1}_{Q(J^{+})}}{\nu(J^+)}+\frac{\mathbf{1}_{Q(J^{-})}}{\nu(J^-)}\right),
\end{equation}
for \[\widetilde{\nu}(J):= \sqrt{\frac{\nu(J^+)\nu(J^-)}{\nu(J^+)+\nu(J^-)}}.\]

For weights $\mu, \sigma,$ define the operator $U_{\mu}:L^2(\mathcal{H},\sigma)\to L^2(\mathcal{H},\sigma)$ by $U_\mu(f):=U(\mu \,f),$ for $f\in L^2(\mathcal{H,\sigma}).$ 

Assume that the operator $U\mathcal{M}_{\overline{d}}^{-\frac{1}{2}}U_{\mu}$ is bounded from $L^2_c(\mathcal{H},\mu)$ to $L^2_c(\mathcal{H}, \nu).$

Then the operator is bounded when tested on the indicator functions of the Carleson tiles, $\mathbf{1}_{T(I)},$ for all $I\in \mathcal{D},$ and so
\begin{equation}\label{E:Forwards-testing}
\| U\mathcal{M}_{\overline{d}}^{-\frac{1}{2}} U_{\mu} \mathbf{1}_{T(I)}\|^2_{L^2_c(\mathcal{H},\nu)} \lesssim \|\mathbf{1}_{T(I)}\|^2_{L^2_c(\mathcal{H},\nu)}.
\end{equation}

We calculate that 
\begin{align} \label{E:Upper-on-indicator}
U\mathcal{M}_{\overline{d}}^{-\frac{1}{2}} U_{\mu} \mathbf{1}_{T(I)} &= U\mathcal{M}_{\overline{d}}^{-\frac{1}{2}} U(\mu \mathbf{1}_{T(I)}) \\
&= \mu_I U\mathcal{M}_{\overline{d}}^{-\frac{1}{2}} U(\mathbf{1}_{T(I)}) \notag \\
&= \mu_I U\mathcal{M}_{\overline{d}}^{-\frac{1}{2}} \mathbf{1}_{Q_{\pm}(I)}\notag \\
&= \mu_I \sum_{K\in \mathcal{D}} \langle \widetilde{\mathbf{1}}_{T(K)}, \mathcal{M}_{\overline{d}}^{-\frac{1}{2}}\mathbf{1}_{Q_{\pm}(I)}\rangle_{L^2(\mathcal{H})} \widetilde{\mathbf{1}}_{Q_{\pm}(K)} \notag \\
&= \frac{\sqrt{2}}{2} \mu_I \sum_{K\subset I^{\pm}} \delta^{(I,K)}|K|^{\frac{1}{2}} d_K \widetilde{\mathbf{1}}_{Q_{\pm}(K)}.\notag
\end{align}

For a function $f=\sum_{J\in \mathcal{D}} \langle f, H_J^{\nu}\rangle_{L^2(\mathcal{H,\nu})} H_J^{\nu},$ we define the projection
\[
Q^\nu_I= \sum_{J\subset I} \langle f, H_J^{\nu}\rangle_{L^2_c(\mathcal{H},\nu)} H_J^{\nu}.
\]
Expanding $Q^{\nu}_I U\mathcal{M}_{\overline{d}}^{-\frac{1}{2}} U_{\mu} \mathbf{1}_{T(I)}$ with respect to the basis $\{H_J^{\nu}\}_{J\in\mathcal{D}}$, we calculate the left-hand side of the inequality~\eqref{E:Forwards-testing}:

\begin{align} 
&\|Q^{\nu}_I U\mathcal{M}_{\overline{d}}^{-\frac{1}{2}} U_{\mu} \mathbf{1}_{T(I)}\|^2_{L^2_c(\mathcal{H},\nu)} \notag \\ &= \frac{(\mu_I)^2}{2} \sum_{J\subset I} \left|\langle \sum_{K\subset I^{\pm}} \delta^{(I,K)}|K|^{\frac{1}{2}} d_K \widetilde{\mathbf{1}}_{Q_{\pm}(K)}, H_J^{\nu} \rangle_{L^2_c(\mathcal{H},\nu)}\right|^2 \notag \\
&=\frac{(\mu_I)^2}{2} \sum_{J\subset I} \left| \sum_{K\subset I^{\pm}} \delta^{(I,K)}|K|^{\frac{1}{2}} d_K \langle \widetilde{\mathbf{1}}_{Q_{\pm}(K)}, H_J^{\nu} \rangle_{L^2_c(\mathcal{H},\nu)}\right|^2 \notag \\
&= \frac{(\mu_I)^2}{2}\left| \sum_{K\subset I^{\pm}} \delta^{(I,K)}|K|^{\frac{1}{2}} d_K \langle \widetilde{\mathbf{1}}_{Q_{\pm}(K)}, H_I^{\nu} \rangle_{L^2_c(\mathcal{H},\nu)}\right|^2  \notag \\
&+ \frac{(\mu_I)^2}{2} \sum_{J\subset I^{\pm}} \left| \sum_{K\subset I^{\pm}} \delta^{(I,K)}|K|^{\frac{1}{2}} d_K \langle \widetilde{\mathbf{1}}_{Q_{\pm}(K)}, H_J^{\nu} \rangle_{L^2_c(\mathcal{H},\nu)}\right|^2. \label{E: Forward-testing-split-second}
\end{align}
We compute that
\begin{equation} \label{E:Inner-product-nu-cases}
\langle \widetilde{\mathbf{1}}_{Q_{\pm}(K)}, H_J^{\nu}\rangle_{L^2_c(\mathcal{H},\nu)} = 
\begin{cases}
     0 & J\subset K, \\
      \frac{\delta^{(J,K)}}{|K|}\frac{\widetilde{\nu}(J)}{\nu(J^{\pm})}\left(-\nu(K^{+})+\nu(K^{-})\right) & K\subset J^{\pm}, \\
      \frac{1}{|K|} \widetilde{\nu}(K) & K=J.
      
\end{cases}
\end{equation}
The first term in the sum~\eqref{E: Forward-testing-split-second} now becomes:
\begin{align*}
&\frac{(\mu_I)^2}{2}\left| \sum_{K\subset I^{\pm}} \delta^{(I,K)}|K|^{\frac{1}{2}} d_K \langle \widetilde{\mathbf{1}}_{Q_{\pm}(K)}, H_I^{\nu} \rangle_{L^2_c(\mathcal{H},\nu)}\right|^2 \\
&= \frac{(\mu_I)^2}{2}\left| \sum_{K\subset I^{\pm}} \frac{d_K}{|K|^{\frac{1}{2}}} \frac{\widetilde{\nu}(I)}{\nu(I^{\pm})}(-\nu(K^{+})+\nu(K^{-}))\right|^2.
\end{align*}
We observe that 
\begin{align*}
& \frac{(\mu_I)^2}{2}\left| \sum_{K\subset I^{\pm}} \frac{d_K}{|K|^{\frac{1}{2}}} \frac{\widetilde{\nu}(I)}{\nu(I^{\pm})}(-\nu(K^{+})+\nu(K^{-}))\right|^2 \\
 &\leq \|Q^{\nu}_I U\mathcal{M}_{\overline{d}}^{-\frac{1}{2}} U_{\mu} \mathbf{1}_{T(I)}\|^2_{L^2_c(\mathcal{H},\nu)} \\
&\leq \|U\mathcal{M}_{\overline{d}}^{-\frac{1}{2}} U_{\mu} \mathbf{1}_{T(I)}\|^2_{L^2_c(\mathcal{H},\nu)} \lesssim \|\mathbf{1}_{T(I)}\|^2_{L^2(\mu)} = \mu_I |I|^2.
\end{align*}
Since, by the definition of $\mu,$ $\mu_I = \mu(T(I))=|I|$, for all $I \in \mathcal{D},$ it follows that
\begin{align*}
     \left| \frac{1}{|I|^{\frac{1}{2}}} \sum_{K\subset I^{\pm}} \frac{d_K}{|K|^{\frac{1}{2}}} \frac{\widetilde{\nu}(I)}{\nu(I^{\pm})}(-\nu(K^{+})+\nu(K^{-}))\right| \lesssim 1.
\end{align*}
This yields the condition~\eqref{E:Condition-A}, in Theorem \ref{T:0101}.
We now consider the second term in~\eqref{E: Forward-testing-split-second}.
Since $\langle \mathbf{1}_{Q_{\pm}(K)}, H_J^{\nu}\rangle_{L^2_c(\mathcal{H},\nu)} \neq 0$ for $K\subset J,$ and in the sum below $J\subset I^{\pm}$ and $K\subset I^{\pm},$ we have that 
\begin{align*}
&\frac{\mu_I^2}{2} \sum_{J\subset I^{\pm}} \left| \sum_{K\subset I^{\pm}} \delta^{(I,K)}|K|^{\frac{1}{2}} d_K \langle \widetilde{\mathbf{1}}_{Q_{\pm}(K)}, H_J^{\nu} \rangle_{L^2_c(\mathcal{H},\nu)}\right|^2 \\
&=\frac{\mu_I^2}{2} \sum_{J\subset I^{\pm}} \left| \sum_{K\subset J^{\pm}} \delta^{(I,K)} d_K \frac{\delta^{(J,K)}}{|K|^{\frac{1}{2}}} \frac{\widetilde{\nu}(J)}{\nu(J^{\pm})}(-\nu(K_+)+\nu(K_-))+\delta^{(I,J)}\widetilde{\nu}(J)\frac{d_J}{|J|^{\frac{1}{2}}}\right|^2.
\end{align*}
Let \[
p(J)= \sum_{K\subset J^{\pm}}   \frac{d_K}{|K|^{\frac{1}{2}}} \frac{\delta^{(J,K)}}{\nu(J^{\pm})}(-\nu(K_+)+\nu(K_-)).
\]
Since $J\subset I^{\pm},$ and $K\subset J^{\pm},$ we have that $\delta^{(I,K)}=\delta^{(I,J)}.$ Thus,
\begin{align*}
&\frac{(\mu_I)^2}{2} \sum_{J\subset I^{\pm}} \left| \sum_{K\subset I^{\pm}} \delta^{(I,K)}|K|^{\frac{1}{2}} d_K \langle \widetilde{\mathbf{1}}_{Q_{\pm}(K)}, H_J^{\nu} \rangle_{L^2_c(\mathcal{H},\nu)}\right|^2 \\
&=\frac{(\mu_I)^2}{2} \sum_{J\subset I^{\pm}} (\widetilde{\nu}(J))^2\left| p(J)+\frac{d_J}{|J|^{\frac{1}{2}}}\right|^2.
\end{align*}
As with the first term in the sum~\eqref{E: Forward-testing-split-second}, assuming that the operator is bounded, we get:
\begin{align*}
\left(\frac{1}{|I|} \sum_{J\subset I^{\pm}} (\widetilde{\nu}(J))^2\left| p(J)+\frac{d_J}{|J|^{\frac{1}{2}}}\right|^2 \right)^{\frac{1}{2}} \lesssim 1.
\end{align*} 
Taking the supremum over all $I\in \mathcal{D},$ we obtain the condition~\eqref{E:Condition-B} of Theorem \ref{T:0101}. 

Finally, under the assumption that $U\mathcal{M}_{\overline{d}}^{-\frac{1}{2}}U$ is bounded between $L^2_c(\mathcal{H},\mu)$ and $L^2_c(\mathcal{H},\nu)$, we prove that the third condition holds. 

It follows that $U^*\mathcal{M}_{d}^{-\frac{1}{2}}U^*_{\nu}$ is bounded as an operator between  $L^2_c(\mathcal{H},\nu)$ and $L^2_c(\mathcal{H},\mu).$ 
Therefore, 
\begin{equation}\label{E:Backwards-testing}
\| U^*\mathcal{M}_{d}^{-\frac{1}{2}} U^*_{\mu} \mathbf{1}_{Q(I)}\|^2_{L^2_c(\mathcal{H},\mu)} \lesssim \|\mathbf{1}_{Q(I)}\|^2_{L^2_c(\mathcal{H},\nu)}.
\end{equation}
We expand $\mathbf{1}_{Q(I)} U^*\mathcal{M}_{d}^{-\frac{1}{2}} U^*_{\mu} \mathbf{1}_{Q(I)} $ with respect to the basis $\{\widetilde{\mathbf{1}}_{T(J)}\}_{J\in \mathcal{D}}$ to calculate the left-hand side of the equation~\eqref{E:Backwards-testing}. We get
\begin{align*}
    &\| \mathbf{1}_{Q(I)}U^*\mathcal{M}_{d}^{-\frac{1}{2}} U^*_{\nu} \mathbf{1}_{Q(I)}\|^2_{L^2_c(\mathcal{H},\mu)}  \\
    &=\sum_{J\in \mathcal{D}} \left|\langle U^*\mathcal{M}_{d}^{-\frac{1}{2}} U^*_{\nu} \mathbf{1}_{Q(I)}, \widetilde{\mathbf{1}}_{T(J)} \rangle_{L^2_c(\mathcal{H},\mu)} \right|^2 \mu_J \\
    &= \sum_{J\in \mathcal{D}} \left|\langle \nu \mathbf{1}_{Q(I)}, U\mathcal{M}_{\overline{d}}^{-\frac{1}{2}}U\widetilde{\mathbf{1}}_{T(J)} \rangle_{L^2_c(\mathcal{H},\mu)} \right|^2 \mu_J \\
    &=\sum_{J\in \mathcal{D}} \left|\langle \nu \mathbf{1}_{Q(I)},  \frac{1}{|J|} \sum_{K\subset J^{\pm}} \delta^{(J,K)}|K|^{-\frac{1}{2}} \overline{d_K} \mathbf{1}_{Q_{\pm}(K)} \rangle_{L^2_c(\mathcal{H},\mu)} \right|^2 \mu_J\\
    &= \sum_{J\in \mathcal{D}} \frac{\mu_J}{|J|^2}\left|\sum_{K\subset J^{\pm}} \delta^{(J,K)}|K|^{-\frac{1}{2}} \overline{d_K} \langle \nu \mathbf{1}_{Q(I)},   \mathbf{1}_{Q_{\pm}(K)} \rangle_{L^2_c(\mathcal{H},\mu)} \right|^2  \\
     &= \sum_{J\in \mathcal{D}}\frac{\mu_J}{|J|^2} \left|\sum_{K\subset J^{\pm}} \delta^{(J,K)}|K|^{-\frac{1}{2}} \overline{d_K} (\nu(K^+)-\nu(K^-) )\right|^2  \\
      &= \sum_{J\in \mathcal{D}}\frac{1}{|J|} \left|\sum_{K\subset J^{\pm}} \delta^{(J,K)}|K|^{-\frac{1}{2}} \overline{d_K} (\nu(K^+)-\nu(K^-) )\right|^2  \\
      &\lesssim \|\mathbf{1}_{Q(I)}\|^2_{L^2_c(\mathcal{H},\nu)} = \nu(I).
\end{align*}
Therefore,
\[
\left(\frac{1}{\nu(I)}\sum_{J\in \mathcal{D}}\frac{1}{|J|} \left|\sum_{K\subset J^{\pm}} \delta^{(J,K)}|K|^{-\frac{1}{2}} \overline{d_K} (\nu(K^+)-\nu(K^-) )\right|^2\right)^{\frac{1}{2}} \lesssim 1.
\]

Taking a supremum over all $I\in \mathcal{D},$ we get that the condition~\eqref{E:Condition-C} in Theorem \ref{T:0101} holds, proving the necessity of the conditions $A$, $B$ and $C.$

\subsection{Proof of sufficiency}
We now prove the forward direction. We assume that conditions $A,$ $B,$ and $C$ hold, and prove that the operator $U\mathcal{M}_{\overline{d}}^{-\frac{1}{2}}U_\mu$ is bounded from $L^2_c(\mathcal{H}, \mu)$ to $L^2_c(\mathcal{H}, \mu)$.

Let $f\in L^2(\mathcal{H}, \mu),$ and $g\in L^2(\mathcal{H}, \nu).$ We expand the functions with respect to the weighted bases of $L^2_c(\mathcal{H}, \mu)$ and $L^2_c(\mathcal{H},\nu)$ introduced earlier as $f= \sum_{I\in \mathcal{D}} \widehat{f}_{\mu}(I) h_I^{\mu}$ and $g= \sum_{J\in \mathcal{D}} \widehat{g}_{\nu}(J) H_\nu^{J}.$
Using the previous calculations of $ U\mathcal{M}_{\overline{d}}^{-\frac{1}{2}}U (\mu \mathbf{1}_{T(I)}),$ given by the equation~\eqref{E:Upper-on-indicator}, we have
\begin{align*}
    &\langle U\mathcal{M}_{\overline{d}}^{-\frac{1}{2}}U (\mu f), g \rangle_{L^2_c(\mathcal{H},\nu)} \\
    &= \sum_{I,J\in\mathcal{D}} \widehat{f}_{\mu}(I)\widehat{g}_{\nu}(J)\langle U\mathcal{M}_{\overline{d}}^{-\frac{1}{2}}U_\mu h_I^{\mu}, H_J^{\nu}\rangle_{L^2(\mathcal{H},\nu)} \\
    &= \sqrt{2}\sum_{I,J\in \mathcal{D}} \widehat{f}_{\mu}(I)\widehat{g}_{\nu}(J) \frac{\sqrt{\mu_I}}{|I|}\sum_{\substack{K\subset I^{\pm} \\ K \subset J}} \delta^{(I,K)} d_K |K|^{\frac{1}{2}}\langle \widetilde{\mathbf{1}}_{Q_{\pm}(K)}, H_\nu^{J}\rangle. 
\end{align*}
Since in the innermost sum $K\subset J,$ and $K\subset I^{\pm},$ it follows that $I\cap J \neq \emptyset.$ Therefore,
\begin{align} \label{E:decompositon-into-sums} &\langle U\mathcal{M}_{\overline{d}}^{-\frac{1}{2}}U (\mu f), g \rangle_{L^2_c(\mathcal{H},\nu)} \\
    &= \sqrt{2}\left(\sum_{I=J}+\sum_{J\subset I^{\pm}}+\sum_{I\subset J^{\pm}} \right)\widehat{f}_{\mu}(I)\widehat{g}_{\nu}(J) \frac{\sqrt{\mu_I}}{|I|}\sum_{\substack{K\subset I^{\pm}  \notag \\ K \subset J}} \delta^{(I,K)} d_K |K|^{\frac{1}{2}}\langle \widetilde{\mathbf{1}}_{Q_{\pm}(K)}, H_\nu^{J}\rangle \notag \\
    &:= T_1 + T_2 + T_3, \notag
\end{align}
where $T_1,$ $T_2,$ and $T_3,$ correspond to the three summands in the second line of the equation~\eqref{E:decompositon-into-sums}.
    We analyze each of the terms separately. 
By the Cauchy-Schwarz inequality, we obtain
\begin{align*}
    &|T_1|=\sqrt{2}\left| \sum_{I\in \mathcal{D}} \widehat{f}_{\mu}(I)\widehat{g}_{\nu}(I) \frac{\mu_I}{|I|} \sum_{K\subset J^{\pm}} \frac{d_K}{|K|^{\frac{1}{2}}} \frac{\widetilde{\nu}(I)}{\nu(I^{\pm})}(-\nu(K^+)+\nu(K^{-}))\right| \\
    &\lesssim \left(\sum_{I\in \mathcal{D}} |\widehat{f}_{\mu}(I)|^2\right)^{\frac{1}{2}}\left(\sum_{I\in\mathcal{D}}|\widehat{g}_{\nu}(I)|^2\right)^\frac{1}{2}\sup_{I\in \mathcal{D}}\left|\sum_{K\subset J^{\pm}} \frac{d_K}{|K|^{\frac{1}{2}}} \frac{\widetilde{\nu}(I)}{\nu(I^{\pm})}(-\nu(K^+)+\nu(K^{-}))\right| \\
    & \lesssim A \|f\|_{L^2(\mathcal{H},\mu)} \|g\|_{L^2(\mathcal{H},\nu)}.
\end{align*}
Next,
\begin{align*} &|T_2|=\sqrt{2}\left|\sum_{I\in \mathcal{D}}\sum_{J\subset I^{\pm}} \widehat{f}_{\mu}(I)\widehat{g}_{\nu}(J) \frac{\sqrt{\mu_I}}{|I|}\sum_{\substack{K\subset I^{\pm} \\ K \subset J}} \delta^{(I,K)} d_K |K|^{\frac{1}{2}}\langle \widetilde{\mathbf{1}}_{Q_{\pm}(K)}, H_\nu^{J}\rangle\right| \\
    &=\sqrt{2}\left|\sum_{I\in \mathcal{D}}\sum_{J\subset I^{\pm}} \widehat{f}_{\mu}(I)\widehat{g}_{\nu}(J) \frac{\sqrt{\mu_I}}{|I|}\left(\sum_{\substack{K\subset I^{\pm} \\ K = J}} + \sum_{\substack{K\subset I^{\pm} \\ K \subset J^{\pm}}}\right)\delta^{(I,K)} d_K |K|^{\frac{1}{2}}\langle \widetilde{\mathbf{1}}_{Q_{\pm}(K)}, H_\nu^{J}\rangle\right| \\
    &= \sqrt{2}\left|\sum_{I \in \mathcal{D}} \sum_{J\subset I^{\pm}} \widehat{f}_{\mu}(I)\widehat{g}_{\nu}(J) \frac{\sqrt{\mu_I}}{|I|}\left(\widetilde{\nu}(J) \delta^{(I,J)}\left(\frac{d_J}{|J|^{\frac{1}{2}}} 
    + p(J)\right)\right)\right|,
\end{align*}
where \[p(J) = \sum_{K\subset J^{\pm}}   \frac{d_K}{|K|^{\frac{1}{2}}} \frac{\delta^{(J,K)}}{\nu(J^{\pm})}(-\nu(K_+)+\nu(K_-)).\]
The last line follows because, when $J\subset I^{\pm},$ the intersection of the set of all $K\subset I^{\pm}$ and those $K\in \mathcal{D},$ such that $K\subset J^{\pm}$ are precisely $K \subset J^{\pm}.$ In the innermost sum, $\delta^{(I,K)}$ then becomes $\delta^{(I,J)}.$ 

We let $\widetilde{f}_\mu = \sum_{I\in \mathcal{D}} \widehat{f}_{\mu}(I) h_I.$ Then $\|\widetilde{f}_\mu\|_{L^2(\mathbb{R})}=\|f\|_{L^2(\mathcal{H},\mu)} $ and $\widetilde{f}_\mu \in L^2(\mathbb{R}).$

Under the assumption $B$ of Theorem \ref{T:0101}, by Fubini's theorem and the Cauchy-Schwarz inequality, it follows that:
\begin{align*}
    &\left|\sum_{I \in \mathcal{D}} \sum_{J\subset I^{\pm}} \widehat{f}_{\mu}(I)\widehat{g}_{\nu}(J) \frac{\sqrt{\mu_I}}{|I|}\left(\widetilde{\nu}(J) \delta^{(I,J)}\left(\frac{d_J}{|J|^{\frac{1}{2}}} 
    + p(J)\right)\right)\right| \\
    &=\left|\sum_{J \in \mathcal{D}} \sum_{J\subset I^{\pm}} \widehat{f}_{\mu}(I)\widehat{g}_{\nu}(J) \frac{\sqrt{\mu_I}}{|I|}\delta^{(I,J)}\left(\widetilde{\nu}(J) \left(\frac{d_J}{|J|^{\frac{1}{2}}} 
    + p(J)\right)\right)\right| \\
    &=\left|\sum_{J \in \mathcal{D}} \widehat{g}_{\nu}(J) \left(\widetilde{\nu}(J) \left(\frac{d_J}{|J|^{\frac{1}{2}}} 
    + p(J)\right)\right)\sum_{J\subset I^{\pm}} \widehat{f}_{\mu}(I)\frac{\delta^{(I,J)}}{\sqrt{|I|}}\right| \\
    &=\left|\sum_{J \in \mathcal{D}} \widehat{g}_{\nu}(J) \left(\widetilde{\nu}(J) \left(\frac{d_J}{|J|^{\frac{1}{2}}} 
    + p(J)\right)\right)\sum_{J\subsetneq I} \widehat{f}_{\mu}(I)h_I(J)\right| \\
    &=\left|\sum_{J \in \mathcal{D}} \widehat{g}_{\nu}(J) \left(\widetilde{\nu}(J) \left(\frac{d_J}{|J|^{\frac{1}{2}}} 
    + p(J)\right)\right)\langle \widetilde{f}_\mu \rangle_{J}\right| \\
    &\leq\left(\sum_{J \in \mathcal{D}} |\widehat{g}_{\nu}(J)|^2\right)^{\frac{1}{2}} \left( \sum_{J \in \mathcal{D}}  \widetilde{\nu}(J)^2 \left(\frac{d_J}{|J|^{\frac{1}{2}}} 
    + p(J)\right)^2\left(\langle \widetilde{f}_\mu \rangle_{J} \right)^2\right) \\
    & \lesssim B \left(\sum_{J \in \mathcal{D}} |\widehat{g}_{\nu}(J)|^2\right)^{\frac{1}{2}} \left(\sum_{J \in \mathcal{D}} |\widehat{f}_{\mu}(J)|^2\right)^{\frac{1}{2}} = B \|g\|_{L^2(\mathcal{H}, \nu)} \|f\|_{L^2(\mathcal{H}, \mu)}. 
\end{align*}
Hence, $|T_2| \lesssim \|g\|_{L^2(\mathcal{H}, \nu)} \|f\|_{L^2(\mathcal{H}, \mu)}.$
Lastly, we bound the third term of the expression~\eqref{E:decompositon-into-sums} given by:
\begin{align*}
    |T_3|=\sqrt{2} \left|\sum_{I\in \mathcal{D}}\sum_{I\subset J^{\pm}} \widehat{f}_{\mu}(I)\widehat{g}_{\nu}(J) \frac{\sqrt{\mu_I}}{|I|}\sum_{K\subset I^{\pm}} \delta^{(I,K)} d_K |K|^{\frac{1}{2}}\langle \widetilde{\mathbf{1}}_{Q_{\pm}(K)}, H_\nu^{J}\rangle\right|.
\end{align*}
Assume $I\in \mathcal{D},$ $I\subset J^{\pm},$ and $K\subset I^{\pm}.$ Then
\begin{align*}
\langle \widetilde{\mathbf{1}}_{Q_{\pm}(K)}, H_J^{\nu}\rangle &= 
\begin{cases}
      \frac{1}{|K|}\frac{\widetilde{\nu}(J)}{\nu(J^{+})}\left(-\nu(K^{+})+\nu(K^{-})\right) & K\subset J^{+}, \\
       \frac{-1}{|K|}\frac{\widetilde{\nu}(J)}{\nu(J^{-})}\left(-\nu(K^{+})+\nu(K^{-})\right) & K\subset J^{-}
\end{cases} \\
&= \frac{1}{|K|}\left(-\nu(K^{+})+\nu(K^{-})\right) \widetilde{\nu}(J) \begin{cases}
    \frac{1}{\nu(J^{+})} & I \subset J^{+}, \\
    \frac{-1}{\nu(J^{-})} & I \subset J^{-}.
\end{cases}
\end{align*}
Let \[
q(I)=\sum_{K\subset I^{\pm}} \delta^{(I,K)} \frac{d_K}{|K|^{\frac{1}{2}}} (-\nu(K^{+})+\nu(K^-)).
\]
Then, 
\begin{align*}
    &\sqrt{2} \left|\sum_{I\in \mathcal{D}}\sum_{I\subset J^{\pm}} \widehat{f}_{\mu}(I)\widehat{g}_{\nu}(J) \frac{\sqrt{\mu_I}}{|I|}\sum_{K\subset I^{\pm}} \delta^{(I,K)} d_K |K|^{\frac{1}{2}}\langle \widetilde{\mathbf{1}}_{Q_{\pm}(K)}, H_\nu^{J}\rangle\right|\\
    &= \sqrt{2} \left|\sum_{I\in \mathcal{D}}\sum_{I\subset J^{\pm}} \widehat{f}_{\mu}(I)\widehat{g}_{\nu}(J) \frac{\sqrt{\mu_I}}{|I|}\frac{\widetilde{\nu}(J)}{\nu(J^{\pm})} \delta^{(J,I)} q(I) \right| \\
    &= \sqrt{2} \left|\sum_{I\in \mathcal{D}} \widehat{f}_{\mu}(I) \frac{\sqrt{\mu_I}}{|I|} q(I) \sum_{I\subset J^{\pm}} \widehat{g}_{\nu}(J) \frac{\widetilde{\nu}(J)}{\nu(J^{\pm})} \delta^{(J,I)} \right| \\
     &= \sqrt{2} \left|\sum_{I\in \mathcal{D}} \widehat{f}_{\mu}(I) \frac{\sqrt{\mu_I}}{|I|} q(I) \left\langle g, \frac{\mathbf{1}_{Q(I)}}{\nu(I)} \right\rangle_{L^2(\mathcal{H},\nu)}\right| \\
    &\leq \sqrt{2}\left(\sum_{I\in \mathcal{D}} |\widehat{f}_{\mu}(I)|^2\right)^{\frac{1}{2}} \left(\sum_{I\in \mathcal{D}} \frac{1}{|I|} (q(I))^2 \left(\left\langle g, \frac{\mathbf{1}_{Q(I)}}{\nu(I)} \right\rangle_ {L^2(\mathcal{H},\nu)}\right)^2 \right)^{\frac{1}{2}}.
\end{align*}
Condition $C$ from Theorem \ref{T:0101} now implies that coefficients $\frac{(q(I))^2}{|I|}$ are Carleson coefficients, and so
\begin{align*}
&|T_3| = \sqrt{2}\left(\sum_{I\in \mathcal{D}} |\widehat{f}_{\mu}(I)|^2\right)^{\frac{1}{2}} \left(\sum_{I\in \mathcal{D}} \frac{1}{|I|} (q(I))^2 \left(\left\langle g, \frac{\mathbf{1}_{Q(I)}}{\nu(I)} \right\rangle_{L^2(\mathcal{H},\nu)}\right)^2 \right)^{\frac{1}{2}} \\ 
&\lesssim \|f\|_{L^2(\mathcal{H}, \mu)} \|g\|_{L^2(\mathcal{H}, \nu)}.
\end{align*}
\subsection{Control of the bounds by the $BMO$ norms}
We note that $\|P_{b,d}^{(0,1,0,1)}\|_{L^2(\mathbb{R})} \approx (A+B+C).$ It follows that $(A+B+C) \lesssim \|b\|_{BMO} \|d\|_{BMO}.$ We now give a direct proof of this fact. The argument for each of the constants relies on $H^1-BMO$ duality.

Suppose that $b,d \in BMO.$ 
 Theorem 1.3 in~\cite{PotCarSawWic} gives conditions for the boundedness of the operator $P_{b,\epsilon}^{(0,1,0,0)}.$ Taking $\epsilon_I= 1,$ for all $I\in \mathcal{D},$ the condition is implied by $b\in BMO.$ We prove this fact directly. 
In particular, we prove that if $I\in \mathcal{D},$ it follows that 
\begin{equation} \label{E:(0,1,1,0)-cond-two-duality}
\left(\sum_{J\subset I}\frac{1}{|J|}\left( \sum_{K\subset J^+}|b_K|^2-\sum_{K\subset J^-}|b_K|^2\right)^2\right)^{\frac{1}{2}} \leq \|b\|_{BMO} \left(\sum_{J\subset I} |b_J|^2\right)^{\frac{1}{2}}.
\end{equation}
Using $l^2$ duality, the equation~\eqref{E:(0,1,1,0)-cond-two-duality} is equivalent to proving that 
\begin{equation*}
\left|\sum_{J\subset I}\frac{1}{\sqrt{|J|}} \left(\sum_{K\subset J^+}|b_K|^2-\sum_{K\subset J^-}|b_K|^2\right)\overline{c_J} \right| \leq \|b\|_{BMO} \left(\sum_{J\subset I} |b_J|^2\right)^{\frac{1}{2}} \left(\sum_{J\subset I}|c_J|^2\right)^{\frac{1}{2}},
\end{equation*}
for a sequence $c^I=\{c_J\}_{J\subset I} $ such that $\|c^I\|_{l^2}< \infty.$
We also use $c^I$ to denote the function $c^I:= \sum_{J\subset I} c_J h_J.$ We use $b^{I^{\pm}}$ to denote the function $b^{I^{\pm}} = \sum_{J\subset I^{\pm}} b_J h_J.$

Now,
\begin{align*}
\left|\sum_{J\subset I}\frac{1}{\sqrt{|J|}} \left(\sum_{K\subset J^+}|b_K|^2-\sum_{K\subset J^-}|b_K|^2\right)\overline{c_J} \right|& = \left|\sum_{J\subset I}\sum_{K\subsetneq J}\frac{\delta^{(J,K)}}{\sqrt{|J|}} |b_K|^2 \overline{c_J}  \right| \\
& = \left|\sum_{K\subsetneq I} |b_K|^2  \sum_{K\subsetneq J\subset I}\frac{\delta^{(J,K)}}{\sqrt{|J|}} \overline{c_J}  \right| \\
& = \left|\sum_{K\subsetneq I} |b_K|^2 \langle \overline{c^I} \rangle_{K} \right| \\
& = \left|\sum_{K\subsetneq I} b_K \overline{b_K} \langle \overline{c^I} \rangle_{K} \right|\\
&=\left|\left \langle b, \sum_{K\subsetneq I} b_K \langle c^I \rangle_{K}\right \rangle_{L^2(\mathbb{R})} \right| \\
& \lesssim \|b\|_{BMO} \|S_{\Phi}\|_{L^1},
\end{align*}
where the last inequality follows by $H^1-BMO$ duality. 
Here,
$\Phi := \sum_{K\subsetneq I} b_K \langle c^I \rangle_K$.
Note that $S\Phi$ denotes the square function defined by:

$S\Phi:= \left(\sum_{I\in \mathcal{D}} (|\Phi_I|)^2 \frac{\mathbf{1}_I}{|I|} \right)^{\frac{1}{2}}$.
Note that 
\begin{equation} \label{E:Square-function-reasoning}
(S\Phi)^2= \sum_{K\subsetneq I} (|b_K|)^2 |\langle c^I \rangle_K|^2\frac{\mathbf{1}_K}{|K|} \leq (M{c^I})^2 (Sb^{I^{\pm}})^2, \end{equation}
where $M$ denotes the Hardy-Littlewood maximal function.
Since the maximal and square functions are bounded on $L^2(\mathbb{R})$, it follows that
\begin{align*}
    \|S\Phi\|_{L^1} \leq \|Mc^I\|_{L^2(\mathbb{R})} \|Sb_{I^{\pm}}\|_{L^2(\mathbb{R})} &\lesssim \|c^{I}\|_{L^2(\mathbb{R})} \|b^{I^{\pm}}\|_{L^2(\mathbb{R})}  \\
    &=\left(\sum_{J\subset I} |c_J|^2\right)^{\frac{1}{2}} \left(\sum_{J\subsetneq I} |b_J|^2\right)^{\frac{1}{2}} \\
    &\leq \left(\sum_{J\subset I} |c_J|^2\right)^{\frac{1}{2}} \left(\sum_{J\subset I} |b_J|^2\right)^{\frac{1}{2}}.
\end{align*}
Thus, the equation~\eqref{E:(0,1,1,0)-cond-two-duality} holds. We now turn to proving that $A \lesssim \|b\|_{BMO} \|d\|_{BMO}.$ We prove that for all $I\in \mathcal{D},$
\begin{equation} \label{E:l^infty-term}
    \frac{1}{\sqrt{|I|}} \left|\sum_{K\subset I^{\pm}} \frac{d_K}{|K^{\frac{1}{2}}|}\frac{\widetilde{\nu}(I)}{\nu(I^{\pm})} (-\nu(K^+)+\nu(K^-))\right| \lesssim \|b\|_{BMO} \|d\|_{BMO}.
\end{equation}
We estimate the terms corresponding to $I^{+}$ and $I^{-}$ in the sum on the left-hand side of the equation~\eqref{E:l^infty-term} separately, and use triangle inequality to obtain the final bound as follows:
\begin{align*}
&\frac{1}{\sqrt{|I|}} \left|\sum_{K\subset I^{+}} \frac{d_K}{|K^{\frac{1}{2}}|}\frac{\widetilde{\nu}(I)}{\nu(I^{\pm})} (-\nu(K^+)+\nu(K^-))\right| \\
&\leq \left(\sum_{J\subset I^{+}}\frac{1}{|K|}(-\nu(K^{+})+\nu(K^{-}))^2\right)^{\frac{1}{2}} \frac{1}{\sqrt{|I|}}\left(\sum_{K\subset I^{+}} |d_K|^2 \frac{(\widetilde{\nu}(I))^2}{(\nu(I^{+}))^2}\right)^{\frac{1}{2}}.
\end{align*}
By equation~\eqref{E:(0,1,1,0)-cond-two-duality}, we have that 
\[\left(\sum_{J\subset I^{+}}\frac{1}{|K|}(-\nu(K^{+})+\nu(K^{-}))^2\right)^{\frac{1}{2}} \lesssim \|b\|_{BMO} (\nu(I^+))^{\frac{1}{2}}. \]
Thus, 
\begin{align*}
&\frac{1}{\sqrt{|I|}} \left|\sum_{K\subset I^{+}} \frac{d_K}{|K^{\frac{1}{2}}|}\frac{\widetilde{\nu}(I)}{\nu(I^{\pm})} (-\nu(K^+)+\nu(K^-))\right| \\
&\lesssim \|b\|_{BMO} \frac{\widetilde{\nu}(I)}{\nu(I^+)}\left(\nu(I^+)\right)^{\frac{1}{2}} \left(\frac{1}{|I|}\sum_{K\subset I^{+}} |d_K|^2 \right)^{\frac{1}{2}} \\
& \leq \|b\|_{BMO} \|d\|_{BMO} \left(\frac{\nu(I^+)\nu(I^-)}{\nu(I)\nu(I^+)}\right)^{\frac{1}{2}} \\
& \leq \|b\|_{BMO}\|d\|_{BMO}.
\end{align*}
Thus, $A\lesssim \|b\|_{BMO}\|d\|_{BMO}. $ The proof that $B+C \lesssim \|b\|_{BMO}\|d\|_{BMO}$ follows similar reasoning and the $H^1-BMO$ duality, so we leave it to the interested reader.\end{proof}

\bibliographystyle{amsplain}

\begin{bibdiv}
\begin{biblist}
\bib{Cal}{article}{
   author={Calder\'{o}n, A. P.},
   title={Commutators of singular integral operators},
   journal={Proc. Nat. Acad. Sci. USA},
   volume={53},
   date={1965},
   pages={1092--1099}
}
\bib{Chu}{thesis}{
   author={Chung, Dae-Won},
   title={Commutators and dyadic paraproducts on weighted Lebesgue spaces},
   school={University of New Mexico},
   date={2010},
   type={Ph.D. dissertation}
}
\bib{Cru}{article}{
   author={Cruz-Uribe, D.},
   title={The product of unbounded Toeplitz operators},
   journal={IEOT},
   volume={20},
   date={1994},
   pages={231--237}
}

\bib{BenMalNai}{article}{
   author={Bényi, Á.},
   author={Maldonado, D.},
   author={Naibo, V.},
   title={What is a paraproduct?},
   journal={Notices Amer. Math. Soc.},
   volume={57},
   date={2010},
   number={7},
   pages={858--860}
}
\bib{KhaNik}{book}{
   editor={Khavin, V. and Nikolskii, N.},
   title={Linear and complex analysis. Problem book 3. Part I},
   series={Lecture Notes in Math.},
   volume={1573},
   publisher={Springer-Verlag},
   address={Berlin},
   date={1994}
}


\bib{Lac}{article}{
   author={Lacey, M.},
   title={Two-weight inequality for the Hilbert transform: A real variable characterization, II},
   journal={Duke Math. J.},
   volume={163},
   date={2014},
   number={15},
   pages={2821--2840}
}

\bib{LacSawSheUri}{article}{
   author={Lacey, M.},
   author={Sawyer, E.},
   author={Shen, C.}, 
   author={Uriarte-Tuero, I.},
   title={Two-weight inequality for the Hilbert transform: A real variable characterization, I},
   journal={Duke Math. J.},
   volume={163},
   date={2014},
   number={15},
   pages={2795--2820}
}

\bib{LacSawUri}{article}{
   author={Lacey, M. T.},
   author={Sawyer, E. T.},
   author={Uriarte-Tuero, I.},
   title={Two weight inequality for the Hilbert transform: a real variable characterization, II},
   journal={Duke Math. J.},
   volume={163},
   date={2014},
   number={15},
   pages={2821--2840}
}
\bib{PotCarSawWic}{article}{
    author={Pott,S.},
    author={Reguera, M.},
    author={Sawyer, E.},
    author={Wick, B.},
    title={Composition of dyadic paraproducts}, 
    journal={Adv. Math.},
    volume={298}, 
    date={2017},
    pages={581-611}
}
\bib{Sar}{article}{
   author={Sarason, D.},
   title={Exposed points in $H^1$, I},
   book={
      series={Oper. Theory Adv. Appl.},
      volume={41},
      date={1989},
   },
   pages={485--496}
}
\bib{Naz}{misc}{
    author={Nazarov, F.},
    title = {A counterexample to Sarason’s conjecture}, 
    year = {1997}, 
    note = {available at \url{http://www.math.msu.edu/~fedja/Preprints/Sarason.ps}},
    pages = {1--17}
}
\bib{TreVol}{article}{
   author={Treil, S.},
   author={Volberg, A.},
   title={Wavelets and the angle between past and future},
   series={J. Funct. Anal.},
   volume={143},
   date={1997},
   number={2},
   pages={269--308}
}

\bib{Wic}{article}{
   author={Wick, B. D.},
   title={Commutators, BMO, Hardy spaces and factorization: A survey},
   journal={Real Anal. Exchange},
   volume={45},
   date={2020},
   number={1},
   pages={1--28}
}


\end{biblist}
\end{bibdiv}

\end{document}